# Estimation of samples relevance by their histograms.

## Mikhail A. Antonets

### 1. Introduction.

Main object, discussed in the paper is a set $R$ of the samples which are the functions on a finite set $T = \{1, 2, \dots, |T|\}$ with values from a finite set $V$. We suppose, that these samples are the result of multiple registration of a dynamic system characteristics that evolve in a mode of operation. In the case, when the function $f$ on the set $T$ was got in a such manner but under conditions, that the mode of operation was not established exactly the following question may be formulated: is the function $f$ relevant to the mode of operation presented by the set of samples $R$.

We suggest an answer to this question using the defined below set $M$ of histograms :

for any function $r$ from the set $R$ its histogram $m$ is the function on the set $V$ defined by the relation

$$m(v) = |\{t: t \in T, r(s) = v\}|, v \in V$$

Any histogram $m$ from the set satisfies $M$ the equality

$$\sum_{v \in V} m(v) = |T| \qquad (1)$$

Our solutions is based on the supporting and covering weights for the set $M$ introduced in [1].

### 2. Supporting weights and covering weights.

**Definition 1.** A weight $x$ on the set $V$ is a nonnegative function such that

$$|x| = \sum_{v \in V} x(v) = 1 \qquad (2)$$

The set of all weights on the set V will be denoted by $\Delta_V$.

For any functions $x, m$ defined on the set $V$ we set

$$(x, m)_V = \sum_{v \in V} x(v) m(v)$$

**Definition 2.** A weight $\underline{x}^M$ on the set $V$ is supporting weight for the set $M$ of the function on the set $V$ if for any weight $x$ on the set $V$ the following inequalities hold

$$min_{m \in M} (\underline{x}^M, m)_V \geq min_{m \in M} (x, m)_V$$

The set of supporting weights $\{\underline{x}^M\}$ coincides with the set of solutions of the variational problem

$$min_{m \in M} (\underline{x}^M, m)_V = max_{x \in \Delta_V} min_{m \in M} (x, m)_V \qquad (3)$$

As a quantitative estimation of the degree of relevance, i.e. the correspondence of the tested function $f$ to the criteria for the formation of the set of functions $R$ is based on, we propose the value

$$\underline{r}^M(f) = \sum_{v \in V} m_f(v) \underline{x}^M(v) \qquad (4)$$

where $\underline{x}^M$ is the constructed weight and $m_f$ is histogram of the function $f$.

**Definition 2.** A weight $\bar{x}^M$ on a set $V$ will be called covering one for a set of functions $M$ if for any weight $x$ on the set $V$ it satisfies the inequality

$$max_{m \in M}(\bar{x}^M, m)_V \leq max_{m \in M}(x, m)_V$$

The set of covering weights $\{\bar{x}^M\}$ coincides with the set of solutions of the variational problem

$$max_{m \in M}(\bar{x}^M, m)_V = min_{x \in \Delta_V} max_{m \in M}(x, m)_V$$

The weight $\bar{x}^M$ characterizes the irrelevance of the function $f$ to the set of functions $R$ by means of the value

$$\bar{s}^M(f) = \sum_{v \in V} m_f(v) \bar{x}^M(v) \qquad (5)$$

which decreases with the relevance's growth.

The below-formulated well-known Theorem 1 and Theorem 2 (see [2,3]) states that the discussed variational problems can be reduced to the following linear programming problems:

**Theorem 1.** Any supporting weight $\underline{x}^M$ is a solution of the following problem:

find all the pairs $\{\underline{\alpha}^M, \underline{x}^M\}$, $\underline{\alpha}^M \in \mathbb{R}$, $\underline{x}^M \in \Delta_V$, that maximize the value $\alpha$ under the conditions (2), the inequalities

$$x(v) \geq 0, v \in V, \qquad (6)$$

and the inequalities

$$\alpha - (x, m)_V \leq 0, \quad \forall m \in M \qquad (7)$$

**Theorem 2.** Any covering weight $\bar{x}^M$ is a solution of the following problem:

Find all the pairs $\{\bar{\alpha}^M, \bar{x}^M\}$, $\bar{\alpha}^M \in \mathbb{R}$, $\bar{x}^M \in \Delta_V$, that minimize the value $\alpha$ under the conditions (2), the inequalities

$$\alpha - (x, m)_V \geq 0, \quad \forall m \in M \qquad (8)$$

Since for the uniformly distributed weight $x_0(v) \equiv \frac{1}{|V|}$ on the set $V$ the following equalities holds

$$\sum_{v \in V} m(v) x_0(v) = \frac{|T|}{|V|}$$

then the following inequality hold

$$\underline{\alpha}^M \geq \frac{|T|}{|V|} \tag{9}$$

$$\overline{\alpha}^M \leq \frac{|T|}{|V|} \tag{10}$$

## 3. Reduction of dimension of variation problems

**Theorem 3.** For an element $w$ from the set $V$ and an arbitrary finite set $M$ of nonnegative functions on the set $V$ the following statements hold:

1) if for any function $m$ from the set $M$ the following inequality take place

$$m(w) < \frac{1}{|V|-1} \sum_{v \in V, v \neq w} m(v) \tag{11}$$

then for any supporting weight $\underline{x}^M$ the following equality holds

$$\underline{x}^M(w) = 0$$

2) if for any function $m$ from the set $M$ the following inequality take place

$$m(w) > \frac{1}{|V|-1} \sum_{v \in V, v \neq w} m(v) \tag{12}$$

then for any covering weight $\overline{x}^M$ the following equality holds

$$\overline{x}^M(w) = 0$$

**Proof.** To prove the assertion 1) for given supporting weight $\underline{x}^M$ let us construct the weight $\tilde{x}$, assuming that

$$\tilde{x}(w) = 0$$

and

$$\tilde{x}(v) = \underline{x}^M(v) + \frac{\underline{x}^M(w)}{|V|-1}$$

for any $v, v \neq w$.

Then

$$\sum_{v \in V, v \neq w} \tilde{x}(v) = 1$$

and for arbitrary function $m$ from the set $M$ the following relations hold

$$(\tilde{x}, m)_V - (\underline{x}^M, m)_V = \sum_{v \in V, v \neq w} m(v)(\underline{x}^M(v) + \frac{\underline{x}^M(w)}{|V|-1}) - \sum_{v \in V} \underline{x}^M(v) m(v) =$$

$$= \underline{x}^M(w)(\sum_{v \in V, v \neq w} m(v) \frac{1}{|V|-1} - m(w))$$

and by virtue inequality (11) the inequality

$$(\tilde{x}, m)_V - (\underline{x}^M, m)_V > 0 \qquad (13)$$

takes place when

$$\underline{x}^M(w) \neq 0$$

in contradiction with supporting weight definition. The proof of the assertion 2) is similarly. The theorem is proved.

**Corollary.** In the case when the set $M$ is a set of histogram, generated by functions on the set $T$ with the values in the set $V$ there is the equality

$$\frac{1}{|V|-1} \sum_{v \in V, v \neq w} m(v) - m(w) = \frac{1}{|V|-1}(|T| - m(w)) - m(w) = \frac{1}{|V|-1}(|T| - |V|m(w))$$

and inequality (11) takes the form

$$m(w) < \frac{|T|}{|V|} \qquad (13)$$

and inequality (12) takes the form

$$m(w) > \frac{|T|}{|V|} \qquad (14)$$

**4. Classification of the set of histograms for $V = \{0, 1\}$.**

The sets $R$ of samples with values from the set $\{0,1\}$ gives the most simple examples of the histogram. Any histogram in this case is the pair of natural number $\{m(0), m(1)\}$ satisfying the equality

$$m(0) + m(1) = |T| \qquad (15)$$

The weights $\underline{x}^M, \overline{x}^M$ and its dual weights $\overline{x}^{\hat{V}}$, $\underline{x}^{\hat{V}}$, $\hat{V} = \{\hat{v} : v \in V\}$, where

$$\hat{v} \colon M \to V, v(m) = m(v)$$

(see [1]) may be constructed in the explicit form.

**Theorem 4.** For any set of histogram $M$ generated by the set $R$ of all functions on the finite set $T$ and the set of value $V = \{0, 1\}$ the following assertions are:

1) Let for any histogram $m$ from the set $M$ the following inequality hold

$$m(0) > m(1) \qquad (16)$$

then the following equality take place.

$$\underline{\alpha}^M = min_{m \in M} m(0) \qquad (17)$$

$$\underline{x}^M = \{1,0\} \tag{18}$$

$$\overline{\alpha}^M = max_{m \in M} m(1) \tag{19}$$

$$\overline{x}^M = \{0,1\} \tag{20}$$

2) if the set $M$ contains such histograms $m', m''$ that the following inequality take place

$$m'(1) \geq m'(0) \tag{21}$$
$$m''(1) \leq m''(0) \tag{22}$$

then the following equalities hold

$$\underline{\alpha}^M = \overline{\alpha}^M = \frac{|T|}{2} \tag{23}$$

$$\underline{x}^M = \overline{x}^M = \{\tfrac{1}{2}, \tfrac{1}{2}\} \tag{24}$$

**Proof.** In the case 1) it follow from the condition (16) that to maximize the sum

$$m(0)x(0) + m(1)x(1)$$

for any function $m$ from the set $M$ it is necessary put $x(1) = 0$ that lead us to equality (17). In the case 2) by virtue inequalities (21), (22) it follows that the value $(x, m')_V$ no decrease and the value $(x, m'')_V$ no increase when the values $x(1)$ grows. Therefore the equalities (15) and

$$\left(\{\tfrac{1}{2}, \tfrac{1}{2}\}, m'\right) = \left(\{\tfrac{1}{2}, \tfrac{1}{2}\}, m''\right) = \frac{|T|}{2}$$

lead to relations (23), (24). The theorem is proved.

The principle of complementary slackness for the weight $\overline{x}^{\hat{V}}$ gives in the case 1) the relation

$$\sum_{m \in M} m(0) \overline{x}^{\hat{V}}(m) = min_{m \in M} m(0) \tag{25}$$

That implies that the set of covering weights for the set of functions $\hat{V}$ consist of all weights vanishing outside the set of all functions from the set $M$ that reach the value $min_{m \in M} m(0)$.

Similarly we get that the set of supporting weights for the set of functions $\hat{V}$ consist of all weights vanishing outside the set of all functions from the set $M$ that reach the value $max_{m \in M} m(1)$.

In the case 2) the principle of complementary slackness for the weight $\overline{x}^{\hat{V}}$ gives the relation

$$\sum_{m \in M} m(j) \overline{x}^{\hat{V}}(m) = \frac{|T|}{2} \quad j = 0,1 \tag{25}$$

For any pair functions $u = \{m', m''\}, m', m'' \in M$ satisfying the condition (21), (22), there exists a solution $\overline{x}^{\hat{V}}_u$ of the system of two linear equations (25), belonging to the set of covering weight for the set $\hat{V}$:

$$\bar{x}^{\widehat{V}}{}_u(m) = 0 \text{ for all m distinct from } m', m'' \qquad (26)$$

$$\bar{x}^{\widehat{V}}{}_u(m') = \frac{m''(0) - \frac{1}{2}|T|}{m''(0) - m'(0)} = \frac{\frac{1}{2}|T| - m''(1)}{m'(1) - m''(1)} \qquad (27)$$

$$\bar{x}^{\widehat{V}}{}_u(m'') = \frac{\frac{1}{2}|T| - m'(0)}{m''(0) - m'(0)} = \frac{m'(1) - \frac{1}{2}|T|}{m'(1) - m''(1)} \qquad (28)$$

**References.**